\documentclass[11pt]{amsart}

\usepackage{amsmath,amsthm, amscd, amssymb, amsfonts}

\usepackage{epsfig}

\newcommand{\sou}{\mathfrak s}
\newcommand{\tgt}{\mathfrak e}
\DeclareMathOperator*{\Tim}{\times}
\newcommand{\ftimes}{\sideset{_\tgt}{_\sou}\Tim}

\def\beqn{\begin{eqnarray}}
\def\eeqn{\end{eqnarray}}
\def\beqn*{\begin{eqnarray*}}
\def\eeqn*{\end{eqnarray*}}

\newcommand{\Vc}{{\mathcal V}}
\newcommand{\Hc}{{\mathcal H}}

\newcommand{\G}{{\mathcal G}}
\newcommand{\D}{{\mathcal D}}
\newcommand{\Ec}{{\mathcal E}}
\newcommand{\Pc}{{\mathcal P}}

\newcommand{\cent}{{\mathcal Z}}

\newcommand{\fde}{{\rightharpoonup}}
\newcommand{\gde}{{\leftharpoonup}}

\newcommand{\ot}{{\otimes}}

\newcommand{\lv}{\lambda_V}
\newcommand{\lh}{\lambda_H}
\newcommand{\rv}{\rho_V}
\newcommand{\rh}{\rho_H}
\newcommand{\qui}{\mbox{\rm Quiv\,}(\mathcal{P})}

\newcommand{\Do}{\Vc\bowtie\Hc}

\newcommand{\s}{{\sigma}}

\newcommand\id{\operatorname{id}}

\theoremstyle{plain}

\numberwithin{equation}{section}

\newtheorem{teo}{Theorem}[section]
\newtheorem{lema}[teo]{Lemma}

\newtheorem{cor}[teo]{Corollary}
\newtheorem{prop}[teo]{Proposition}

\theoremstyle{definition}

\newtheorem{defi}[teo]{Definition}

\theoremstyle{remark}

\newtheorem{rmk}[teo]{Remark}

\def\pf{\begin{proof}}

\def\epf{\end{proof}}

\theoremstyle{remark}

\begin{document}

\title[On Braided groupoids]{On braided groupoids}
\author[Maldonado  and Mombelli]{Carolina Maldonado and Juan Mart\'\i n Mombelli
}
\thanks{This work was partially supported by
Agencia C\'ordoba Ciencia, ANPCyT-Foncyt, CONICET, Fundaci\'on
Antorchas and Secyt (UNC)}
\address{Facultad de Matem\'atica, Astronom\'\i a y F\'\i sica
\newline \indent
Universidad Nacional de C\'ordoba
\newline
\indent CIEM -- CONICET
\newline \indent Medina Allende s/n
\newline
\indent (5000) Ciudad Universitaria, C\'ordoba, Argentina}
\email{cmaldona@mate.uncor.edu and mombelli@mate.uncor.edu}

\begin{abstract} We study and give examples of braided groupoids,
and {\it a fortiori}, non-degenerate solutions of the
quiver-theoretical braid equation.
\end{abstract}

\date{\today}
\maketitle

\section*{Introduction}

Let $V$ be a vector space over some field and let $R:V\ot V\to
V\ot V$ be a linear operator. One says that $R$ is a solution of
the Quantum Yang-Baxter equation (QYBE, for short) if
$$R^{12}R^{13}R^{23}=R^{23}R^{13}R^{12},$$
where as usual $R^{12} = R\ot \id$, and so on. The study of
solutions of the QYBE, motivated by problems in statistical
mechanics and low dimension topology, has been a central theme in
algebra along the last 25 years. If $R$ is a solution of the QYBE
and $\tau:V\ot V\to V\ot V$ denotes the usual transposition, then
$c:=R\tau$ is a solution of the {\it braid equation}, that is
\begin{equation}\label{braideqn}(c\ot\id)(\id\ot c)
(c\ot\id)=(\id\ot c)(c\ot\id)(\id\ot c). \end{equation}

Thus, there is a bijective correspondence between solutions of the
QYBE and solutions of the braid equation.

\bigbreak Drinfeld observed in \cite{D} that both the QYBE and the
braid equation have sense if $V$ is just a set and $R: V\times V
\to V\times V$ is just a map; again, there is a bijective
correspondence between solutions of one and the other. He called
this the set-theoretical QYBE and proposed its study as a
meaningful problem. Note that any solution of the set-theoretical
QYBE gives rise, by linearization, to a solution of the QYBE in
the category of vector spaces. Drinfeld's problem was attacked by
two groups of mathematicians: Etingof-Schedler-Soloviev, see
\cite{ESS, S}, and Lu-Yan-Zhu, see \cite{LYZ, lyz3}. See also
\cite{EGS}, where indecomposable solutions on sets with $p$
elements, $p$ a prime, are classified. Later, Takeuchi gave an
alternative presentation of the results by
Etingof-Schedler-Soloviev and Lu-Yan-Zhu, with braided groups
playing a central r\^ole. See \cite{T}.

\bigbreak Now, the braid equation \eqref{braideqn} has sense in
any monoidal category. Another natural monoidal category to
consider is the category $\qui$ of quivers over a fixed set $\Pc$
with tensor product given by pull-back. The braid equation in
$\qui$ is called the  \emph{quiver-theoretical QYBE}, by abuse of
notation. A solution of the braid equation in $\qui$ is called a
\emph{braided quiver}. Note that any finite solution of the
quiver-theoretical QYBE gives rise, by linearization, to a
solution of the QYBE in the category of bimodules over a
commutative separable algebra.

\bigbreak The problem of characterizing solutions of the braid
equation in $\qui$ was attacked by Andruskiewitsch, see \cite{A}.
In particular Theorem 3.10 in {\it loc. cit.} shows that there is
a bijective correspondence between
\begin{itemize}
    \item Non-degenerate braided quivers $\mathcal{A}$,
    \item pairs $(\G, \mathcal{A})$, where $\G$ is a braided
    groupoid and $\mathcal{A}$ is a representation of $\G$ with certain
    properties.
\end{itemize}

In other words, braided groupoids are the fundamental piece of
information in the classification of solutions of the
quiver-theoretical QYBE. This raises naturally the question of
classifying (or at least characterizing) braided groupoids. This
is the problem considered in the present paper.

\bigbreak  Although braided groupoids appear naturally, by the
result quoted above, no systematic investigation of their
structure was undertaken up to now.
 In the paper \cite{AN} a description of matched pair of
groupoids in group-theoretical terms is obtained. See also
\cite[Thm. 3.1]{AM}. The main idea of this work is to use this
result to describe braided groupoids in terms of group theory.

\bigbreak This paper is intended to be as self-contained as
possible. For this reason we include in section \ref{grupoide}
some basics definitions concerning groupoids. In section \ref{s1}
we recall the definition of matched pair of groupoid. We explain
how to obtain matched pairs of groupoids from a collection
$(D,V,H,\gamma)$, where $V,H$ are subgroups of a finite group $D$
such that $V$ intersects trivially any conjugate of $H$ and
$\gamma: V\backslash D/H \to D$ is a section of the canonical
projection. To such collection we attach maps $\lambda_V,
\lambda_H, \rho_V, \rho_H, \triangleright, \triangleleft$
governing the multiplication of $D$, with certain cohomological
flavor.
 In section \ref{br} we recall the definition of braided
groupoid.

\bigbreak Our main result is Theorem \ref{main}, where we
characterize braided groupoids in terms of collections
$(D,V,H,\gamma)$ as before, subject to some restrictions on the
maps $\lambda_V, \lambda_H, \rho_V, \rho_H, \triangleright,
\triangleleft$. In section \ref{exam} we apply the main result to
obtain examples under suitable restrictions. Notably, we analyze
in subsection \ref{handy} a class of braided groupoids that we
call \emph{handy} and give a complete characterization of them in
terms of data including certain "non-associative" group
structures. We stress that such structures appear also in some
other works in the area \cite{N, B}.

In the next subsection, explicit examples of non-handy braided
groupoids are also presented. Finally in section \ref{calculo} we
compute the braiding for the examples given in section \ref{exam}.

\subsection*{Acknowledgment} We are very grateful to Nicol\'as
Andruskiewitsch for his encouragement and comments on a previous
version of this paper. We also wish to thank Sonia Natale for
interesting conversations.

\section{Braided Groupoids}

\subsection{Groupoids}\label{grupoide}

Recall that a (finite) {\it groupoid} is a small category (with
finitely many arrows), such that every morphism has an inverse. We
shall denote a groupoid by $\tgt, \sou:{\mathcal
G}\rightrightarrows {\mathcal P}$, or simply by $\G$, where $\G$
is the set of arrows, $\Pc$ is the set of objects and $\tgt, \sou$
are the target and source maps.

The set of arrows between two objects $P$ and $Q$ is denoted by
$\G(P,Q)$ and we shall also denote $\G(P):=\G(P,P).$ The
composition map is denoted by $m:\G\ftimes{\mathcal G}\rightarrow
\G$, and for two composable arrows $g$ and $h$, that is
$\tgt(g)=\sou(h)$, the composition will be denoted by
juxtaposition: $m(g,h)=gh$.

\medbreak

A {\it morphism} between two groupoids is a functor of the
underlying categories. Two morphisms of groupoids $\phi,\psi:
{\mathcal G}\rightarrow {\mathcal H}$ are {\it similar}, denoted
$\phi\sim\psi$, if there is a natural transformation between them;
that is, if there exists a map $\tau:{\mathcal P}\rightarrow
{\mathcal H}$ such that
$$\phi(g)\tau(\tgt(g))=\tau(\sou(g))\psi(g), \qquad g\in {\mathcal G}.$$

\medbreak

Two groupoids ${\mathcal G}$, ${\mathcal H}$ are {\it isomorphic},
and we write $\G\cong {\mathcal H}$, if there are morphisms $\phi:
{\mathcal G}\rightarrow {\mathcal H}$, $\psi: {\mathcal
H}\rightarrow {\mathcal G}$ such that $\phi\circ\psi$ and
$\psi\circ\phi$ are similar to the corresponding identities.

\medbreak

Any groupoid $\G$ gives rise to a relation on the base $\Pc$, $P
\approx_{\G} Q$ if $\G(P,Q) \neq \emptyset$. A groupoid $\tgt,
\sou :{\mathcal G}\rightrightarrows {\mathcal P}$ is  {\it
connected } if $P\approx_{\G} Q$ for all $P, Q\in{\mathcal P}$.

\medbreak

Let $S$ be an equivalence class in $\Pc$ and let $\G_S$ denote the
corresponding connected groupoid with base $S$; that is,
$\G_S(P,Q)=\G(P,Q)$ for any $P,Q\in S$. Then the groupoid
$\mathcal G$ is isomorphic to the disjoint union of the connected
groupoids $\G_S$: $\G \cong \coprod_{S \in \Pc/\approx} \G_S$.

\medbreak

If $\Hc$ and $\G$ are two isomorphic groupoids over the same base
$\Pc$ then there are (non-canonical) isomorphisms $\G(P)\cong
\Hc(P)$ for all $P\in \Pc$.

\medbreak A subgroupoid $\Hc$ of a groupoid $\G$ is \emph{wide} if
$\Hc$ has the same base $\Pc$ as $\G$.

\medbreak

%Let $\Hc$ and $\G$ be connected isomorphic groupoids and
%$O\in\Pc$. Then the groups $\G(O)$ and $\Hc(O)$ are isomorphic.

Let $\G\rightrightarrows\Pc$ be a groupoid. If $p:\Ec\to \Pc$ is a
map, a {\it left action} of $\G$ to $(\Ec,p)$ is a map $\fde:\G
{\,}_{\tgt}\times_p \Ec \to \Ec$ such that
\begin{equation}\label{ai} p(g\fde x)=\sou(g),\qquad
 g \fde(h \fde x)=gh \fde x,\qquad\id_{p(x)} \, \fde  x =  x,
\end{equation}
for all composable $g,h\in\G$, $x\in\Ec$. Similarly, a {\it right
action} of $\G$ to $(\Ec,p)$ is a map $\gde:
\Ec{\,}_{p}\times_{\sou}\G
 \to \Ec$ such that
\begin{equation}\label{ad} p(x\gde g)=\tgt(g),\qquad
  (x \gde g)\gde h= x\gde gh,\qquad    x\gde \id_{p(x)}   =  x,
\end{equation}
for all composable $g,h\in\G$, $x\in\Ec$

\subsection{Matched Pairs of Groupoids}\label{s1} We briefly recall some
facts about matched pairs of groupoids. See \cite{Ma}, \cite{AA}
and references therein.

\medbreak

A {\it matched pair of groupoids} is a collection
$(\Vc,\Hc,\fde,\gde)$, where  $ \tgt, \sou:\Vc\rightrightarrows
\Pc$ and $\tgt, \sou:{\mathcal H}\rightrightarrows \Pc$ are two
groupoids over the same base $\Pc$, $\rightharpoonup:{\mathcal
H}_{\,\tgt\!}\times_{\sou} \Vc\rightarrow \Vc$ is a left action of
$\Hc$ on $(\Vc, \sou)$, $\leftharpoonup:{\mathcal
H}_{\,\tgt\!}\times_{\sou} \Vc\rightarrow \Hc$ is a right action
of $\Vc$ on $(\Hc, \tgt)$ such that
\begin{equation*} \label{m7}  \tgt(x\rightharpoonup g)=
\sou(x\leftharpoonup g),\quad x\fde gh =(x\fde g)((x\gde g)\fde
h),
\end{equation*}
\begin{equation*}
 xy\gde g = (x\gde (y\fde g))(y\gde g),
\end{equation*}
for composable elements $x,y\in\Hc$ and $g,h\in\Vc$. \medbreak

Let $(\Vc,\Hc,\fde,\gde)$ be a matched pair of groupoids. There is
an associated {\it diagonal groupoid} $\Do$ with set of arrows
$\Vc _{\,\tgt\!}\times_{\sou} \Hc$, base $\Pc$, source, target,
composition and identity given by
\begin{align*} &\sou(g,x)= \sou(g),\quad  \tgt(g,x)= \tgt(x), \\
&(g,x)(h,y)=(g(x\fde h),(x\gde h)y), \quad \mathbf{\id}_P
=(\id_P,\id_P),
\end{align*}
$g,h\in \Vc$, $x,y\in \Hc$, $P\in \Pc$. The groupoids $\Vc$ and
$\Hc$ can be seen as wide subgroupoids of $\Do$. Then we have an
exact factorization of groupoids $\Do = \Vc\Hc$, that is; for
every $z\in \Do $ there are unique elements $x\in \Vc, g\in \Hc$
such that $z= xg$. Conversely, if $\D = \Vc\Hc$ is an exact
factorization of groupoids then there are actions $\gde$, $\fde$
such that $(\Vc,\Hc,\fde,\gde)$ form a matched pair of groupoids,
and $\D \simeq \Do$.

\medbreak

Let us fix a connected groupoid $\D\rightrightarrows\Pc$ and a
point $O\in \Pc$. Set $D=\D(O)$. For each $P\in\Pc$ we fix
$\tau_P\in\D(O,P)$.

\medbreak

In the following we shall study exact factorizations $\D=\Vc\Hc$
where $\Vc$ and $\Hc$ are connected wide subgroupoids. In this
case we can assume that $\tau_P\in \Vc(O,P)$. There is no harm to
assume that $\tau_O=1$. We shall denote $V=\Vc(O), H=\Hc(O)$.

\medbreak

The following lemma will be useful to describe examples of braided
grou\-poids in group-theoretical terms.

\begin{lema}\label{l2} Under the above considerations there
is a bijection between the following data.
\begin{itemize}
\item[i)] Exact factorizations $\D=\Vc\Hc$, where $\Vc, \Hc$ are
connected wide subgroupoids of $\D$,

\item[ii)]  matched pair of groupoids $(\Vc,\Hc,\rightharpoonup,
\leftharpoonup )$ with $\Vc$, $\Hc$ connected, such that $\D\cong
\Vc\bowtie\Hc$ and

\item[iii)]  collections $(V,H,\gamma)$ where $G, H$ are subgroups
of $D$, $\gamma: \Pc\to  D $ is a  (necessarily) injective  map,
and the following conditions are fulfilled
\begin{align} \label{eq1} & D=\coprod_{P\in\Pc} V\gamma_P H,\\
\label{eq2} & V\bigcap zHz^{-1}=\{1\}
\end{align}
for every $z\in D$.
\end{itemize}
\end{lema}

We shall say that the collection $(D,V,H,\gamma)$ satisfying
conditions of Lemma \ref{l2} (iii) is associated to the matched
pair $(\Vc,\Hc,\rightharpoonup, \leftharpoonup )$ or,
equivalently, to the exact factorization $\D=\Vc\Hc$.

\begin{proof} For the implications (i) $\Leftrightarrow$ (ii)
and (ii) $\Rightarrow$ (iii) see \cite[Thm. 3.1]{AM}.

\medbreak

Assume now that $V,H$ are subgroups of $D$ and $\gamma: \Pc\to D$
is a map such that equations \eqref{eq1}, \eqref{eq2} are
fulfilled. Define the wide subgroupoids $\Vc$ and $\Hc$ by
$$\Hc(P,Q):=\tau^{-1}_P\gamma_P H \gamma^{-1}_Q\tau_Q, \quad
\Vc(P,Q):=\tau^{-1}_P V \tau_Q. $$ By construction $\D=\Vc\Hc$ is
an exact factorization.
\end{proof}

\medbreak

\begin{rmk} We can always assume that $\gamma_O=1$.
\end{rmk}

\begin{rmk} Observe that under conditions of Lemma \ref{l2} (iii)
there is a bijection $\Pc\cong V\backslash D/H$ and via this
identification the map $\gamma$ is a section of the canonical
projection. Conditions \eqref{eq1}, \eqref{eq2} imply that $\mid D
\mid= \mid V \mid \mid H \mid \# \Pc$.
\end{rmk}
\medbreak

Summarizing, to obtain an exact factorization of connected
groupoids we need a group $D$, two subgroups $V$ and $H$ of $D$
such that $V$ intersects trivially all conjugates of $H$. Take
$\Pc$ the set of double cosets $V\backslash D/H$ and
$\gamma:\Pc\to D$ is any section of the canonical projection. Some
examples of such collections are the following:
\begin{itemize}
    \item $V, H$ subgroups of $D$ with coprime orders,
    \item  $D=VC$ an exact factorization of groups and $H$ is
    a subgroup of $C$.
\end{itemize}

\medbreak

The following basic observation will be used repeated times.

\begin{lema}\label{unique} Assume that $(D,V,H,\gamma)$ is a
collection satisfying the conditions of Lemma \ref{l2} (iii), then
for any $z\in D$ there exists $g\in V, x\in H$ and $P\in \Pc$
uniquely determined such that $z=g \gamma_P x.$
\end{lema}
\begin{proof} The existence is clear. Assume that $g^{'}
\gamma_Q x^{'}=g \gamma_P x$, then $P=Q$ and $g^{-1}g^{'}=\gamma_P
x x^{'-1} \gamma^{-1}_P\in V\bigcap \gamma_P H\gamma^{-1}_P$,
hence $g=g^{'}$ and $x=x^{'}$.
\end{proof}

\bigbreak

Assume that $(D,V,H,\gamma)$ is associated to the matched pair
$(\Vc,\Hc,\rightharpoonup, \leftharpoonup )$. Thanks to Lemma
\ref{unique} we shall introduce a family of maps. In the next
section these maps will be used to write conditions for a groupoid
to be braided. Concretely, the maps are
$$\triangleright:H\times V\to V, \,\, \triangleleft:H\times V\to
H,$$
$$(\,;\,):H\times V\to \Pc, $$
such that
\begin{equation}\label{b1} xg=(x\triangleright g) \gamma_{(x;g)}
(x\triangleleft g),
\end{equation}
for all $x\in H, g\in V$. Let us also define maps
$$\lv:\Pc\times V\times \Pc\to V, \,\,\, \rv:\Pc\times
V\times \Pc\to H,$$
$$(\,;\,\,\, ;\, ):\Pc\times V\times \Pc\to \Pc$$
and maps
$$\lh:\Pc\times H\times \Pc\to V, \,\,\, \rh:\Pc\times
H\times \Pc\to H, $$
$$<\,;\,\,\, ;\, >:\Pc\times H\times \Pc\to \Pc $$
such that
\begin{equation}\label{b2} \gamma_P g \gamma_Q=\lv(P,g,Q)
\gamma_{(P;g;Q)} \rv(P,g,Q),
\end{equation}
\begin{equation}\label{b3}\gamma_P x \gamma_Q=\lh(P,x,Q)
\gamma_{<P;x;Q>}\rh(P,x,Q),
\end{equation}
for all $P, Q\in\Pc$, $g\in V, x\in H$.

\bigbreak

In the next section we shall study exact factorizations
$\D=\Vc\Hc$ with $\Vc\cong \Hc$. In that case the groups $V$, $H$
are isomorphic.

\medbreak

If $(D,V,H,\gamma)$ is associated to the exact factorization
$\D=\Vc\Hc$, and $\phi: H\to V$ is an isomorphism we shall also
denote by $\phi$ the isomorphism $\phi:\Hc\to \Vc$ given by
$$\phi(\tau^{-1}_P \gamma_P g \gamma^{-1}_Q
\tau_Q)=\tau^{-1}_P\phi(g)\tau_Q .$$

Given such an isomorphism $\phi$, we define the map $m:\D\to\Vc$
as the composition

\begin{align}\label{mult} \D\stackrel{\simeq}{\longrightarrow} \Do
\stackrel{\id\times \phi}{\longrightarrow}
\Vc{\,}_{\tgt}\times_{\sou} \Vc \stackrel{\mu}{\longrightarrow}
\Vc,
\end{align}
where $\mu:\Vc{\,}_{\tgt}\times_{\sou} \Vc
\stackrel{}{\longrightarrow} \Vc$ is the composition.

\medbreak

Using Lemma \ref{l2} the map $m:\D\to \Vc$ can be explicitly
written as follows.

\begin{lema}\label{m} Let $\alpha\in\D(P,Q)$, if
$\alpha=\tau^{-1}_P g\gamma_R x \tau_Q$ for some $g\in V, x\in H,
R\in\Pc$ then
$$m(\alpha)=\tau^{-1}_P g\lh(R,x,Q)\phi(\rh(R,x,Q))\tau_Q.$$
As a particular case if $\alpha\in\D(O,O)$, $\alpha=g\gamma_R x$
then $m(\alpha)=g\phi(x)$.
\end{lema}

\begin{proof} If we have a decomposition $\alpha=\beta_1\beta_2$ where
$\beta_1\in\Vc, \beta_2\in\Hc$ then, by definition,
$m(\alpha)=\beta_1\phi(\beta_2)$. Note that if $\alpha=\tau^{-1}_P
g\gamma_R x \tau_Q$ then
\begin{align*} \alpha&=\tau^{-1}_P
g\gamma_R x\gamma_Q\gamma^{-1}_Q\tau_Q=\tau^{-1}_P g
\lh(R,x,Q)\gamma_{<R;x;Q>}\rh(R,x,Q)\gamma^{-1}_Q\tau_Q\\
&=\tau^{-1}_P
g\lh(R,x,Q)\tau_{<R;x;Q>}\tau^{-1}_{<R;x;Q>}\gamma_{<R;x;Q>}
\rh(R,x,Q)\gamma^{-1}_Q\tau_Q
\end{align*}
where
\begin{align*} \tau^{-1}_P
g\rh(R,x,Q)\tau_{<R;x;Q>}&\in\Vc(P,<R;x;Q>),\\
\tau^{-1}_{<R;x;Q>}\gamma_{<R;x;Q>} \lh(R,x,Q)\gamma^{-1}_Q\tau_Q
&\in\Hc(<R;x;Q>,Q).
\end{align*}
Therefore
\begin{align*} m(\alpha)&=\tau^{-1}_P
g\lh(R,x,Q)\tau_{<R;x;Q>} \phi(\tau^{-1}_{<R;x;Q>}\gamma_{<R;x;Q>}
\rh(R,x,Q)\gamma^{-1}_Q\tau_Q )\\
&=\tau^{-1}_P g\lh(R,x,Q)\phi(\rh(R,x,Q))\tau_Q.
\end{align*}
Since for all $R\in \Pc$, $x\in H$ $\lh(R,x,O)=1,$ and $
\rh(R,x,O)=x$ the second assertion follows.
\end{proof}

\subsection{Braided Groupoids}\label{br}
The notion of {\it braided groupoid} was introduced in \cite{A} in
order to study the quiver-theoretical Yang-Baxter equation.

\begin{defi}[\cite{A}] A {\it braided groupoid} is a collection
$(\Vc, \rightharpoonup, \leftharpoonup)$ where
$\Vc\rightrightarrows\Pc$ is a groupoid,
$(\Vc,\Vc,\rightharpoonup,\leftharpoonup )$ is a matched pair of
groupoids and for every pair $(f,g)\in\Vc_{\tgt}\times_{\sou}\Vc$
the following equation holds:
\begin{equation}\label{ec1}
fg=(f\rightharpoonup g)(f\leftharpoonup g).
\end{equation}
\end{defi}

If $(\Vc, \rightharpoonup, \leftharpoonup)$ is a braided groupoid
then the map $c:\Vc_{\tgt}\times_{\sou}\Vc\to
\Vc_{\tgt}\times_{\sou}\Vc$ defined by
\begin{equation}\label{trenza}c(\alpha,\beta)=(\alpha\rightharpoonup
\beta, \alpha \leftharpoonup \beta)\end{equation} satisfies the
braid equation.

\medbreak

Let $(\Vc,\Hc,\fde,\gde)$ be a matched of groupoids, and
$\phi:\Hc\to \Vc$ a groupoid isomorphism,  recall the diagonal
groupoid $\D$ and the map $m:\D\to \Vc$ as in the previous
section.

\medbreak

Associated to this matched pair of groupoids there is  a new pair
of actions (that we denote with the same symbol) $
\rightharpoonup, \leftharpoonup:\Vc _{\,\tgt\!}\times_{\sou}
\Vc\to \Vc$, and they are defined by
$$g\rightharpoonup h:=\phi^{-1}(g)\rightharpoonup h, \quad g\leftharpoonup
h:=\phi\left( \phi^{-1}(g)\leftharpoonup h\right),$$ for all
composable $g,h\in G$. Since $\phi$ is a groupoid morphism, the
collection $(\Vc,\Vc,\fde,\gde)$ is a matched pair of groupoids.

\begin{lema}\label{l1} The following statements are equivalent.
\begin{itemize}
\item[i)]  $(\Vc, \fde,\gde)$ is a  braided groupoid, \item[ii)]
the map $m:\D \to \Vc$ is a groupoid morphism.
\end{itemize}
\end{lema}
\begin{proof} Let $\mu:\Vc{\,}_{\tgt}\times_{\sou} \Vc \stackrel{}{\longrightarrow}
\Vc$ be the composition. Since $m=\zeta (\id_{\Vc}\ot \phi) \mu$,
where $\zeta:\D\stackrel{\cong}{\longrightarrow} \Do,$ and
$(\id_{\Vc}\ot \phi)$ is a groupoid morphism, then $m$ is a
groupoid morphism if and only if $\mu$ is a groupoid morphism.
Then the proof follows from \cite[Lemma 2.9]{A}, where it is
proven that $(\Vc, \fde,\gde)$ is braided if and only if the
composition map $\mu$ is a groupoid morphism.
\end{proof}

\medbreak

Without lose of generality we can assume that the groupoid $\Vc$
is connected. If $\Vc$ is not connected then $\Vc$ is similar to
the disjoint union of connected groupoids
$$\Vc \cong \coprod_{S \in \Pc/\approx} \Vc_S.$$

\begin{lema} With the notation above $\Vc$ is braided
if and only if  for any $S\in \Pc/\approx$ $\Vc_S$ is a braided
groupoid.
\end{lema}
\begin{proof} The sufficiency is clear. Assume that $\Vc$ is
braided. We only need to show that, for any $S\in \Pc/\approx$,
$\Vc_S$ is stable under the actions $\rightharpoonup,
\leftharpoonup$.

\medbreak

Let $f,g\in\Vc_S$. Using \eqref{ai}, \eqref{ad} we know that
$$\sou(f\rightharpoonup g)=\sou(f),\quad \tgt(f\leftharpoonup g)=\tgt(g).$$
Since $\sou(f), \tgt(g)\in S$ then $f\rightharpoonup g,
f\leftharpoonup g\in\Vc_S$.
\end{proof}

\begin{defi} We shall say that $(D,V,H,\gamma)$ is a
{\it braided groupoid datum} if the associated connected groupoid
$\Vc$ is braided, or, equivalently if the map $m:\D \to \Vc$ is a
groupoid morphism.

\end{defi}

\begin{rmk} The matched pair $(\Vc,\Vc,\fde,\gde)$ and the map
$m:\D\to\Vc$ both depend on the choice of the isomorphism $\phi$.
Sometimes the isomorphism $\phi$ will be clear from the context.
We shall denote $(D,V,H,\gamma, \phi)$ when special emphasis is
needed.
\end{rmk}

\medbreak

The next result gives necessary and sufficient conditions on the
collection $(D,V,H,\gamma,\phi)$ to be a braided groupoid datum.

\begin{teo}\label{main} The collection $(D,V,H,\gamma,\phi)$ is a
braided groupoid datum if and only if
\begin{align}\label{c1} g&=\lv(P,g,Q)\;\phi(\rv(P,g,Q)),\\
\label{c2} \phi(x)&=\lh(P,x,Q)\,\phi(\rh(P,x,Q)),\\
\label{c3} \phi(x)g&=(x\triangleright g)\,\phi(x\triangleleft g),
\end{align}
for all $P,Q\in\Pc$, $g\in V$, $x\in H$.
\end{teo}

\begin{proof} Assume that $(D,V,H,\gamma)$ is a braided groupoid datum. Set
$\alpha=\gamma_Pg\gamma_Q$ $=\lv(P,g,Q)$
$\gamma_{(P;g;Q)}\rv(P,g,Q)$, then using Lema \ref{m} we have that
$m(\alpha)=\lv(P,g,Q)$ $\phi(\rv(P,g,Q))$. Since $m$ is a groupoid
morphism then $m(\alpha)=m(\gamma_P)m(g)$ $m(\gamma_Q)=g$, hence
we have proved equation \eqref{c1}. Equations \eqref{c2},
\eqref{c3} follows in a similar way using equations \eqref{b1},
\eqref{b3}.

\medbreak

Suppose that equations \eqref{c1}, \eqref{c2}, \eqref{c3} are
fulfilled. Let $\alpha$, $\beta \in \D$ two composable elements,
then $\alpha=\tau^{-1}_P g\gamma_R x\tau_Q$,
$\beta=\tau^{-1}_Qh\gamma_S y\tau_M$ for some $g,h\in V$, $x,y\in
H$ and $P,Q,M,R,S\in\Pc$. We shall prove that
$m(\alpha\beta)=m(\alpha)m(\beta)$.

\medbreak

Lemma \ref{m} together with equation \eqref{c2} implies that
$$m(\alpha)=\tau^{-1}_Pg\phi(x)\tau_Q, \quad m(\beta)=\tau^{-1}_Q
h\phi(y)\tau_M.
$$

Let us compute $\alpha\beta$. Define the elements $X,Y\in\Pc$ by
$X:=(R;(x\vartriangleright h);(x;h))$,
$Y:=<X;\rv(R,(x\vartriangleright h),(x;h))(x\vartriangleleft
h);S>$, then

\begin{align*} \alpha\beta&=\tau^{-1}_P g\gamma_R x h\gamma_S y
\tau_M=\tau^{-1}_P g\gamma_R (x\vartriangleright
h)\gamma_{(x;h)}(x\vartriangleleft h)\gamma_S y
\tau_M\\
&=\tau^{-1}_P g\lv(R,(x\vartriangleright
h),(x;h))\gamma_{X}\rv(R,(x\vartriangleright
h),(x;h))(x\vartriangleleft h)\gamma_S y
\tau_M\\
&=\tau^{-1}_P g\lv(R,(x\vartriangleright
h),(x;h))\lh(X,\rv(R,(x\vartriangleright
h),(x;h))(x\vartriangleleft
h),S)\gamma_Y\\
&\rh(X,\rv(R,(x\vartriangleright h),(x;h))(x\vartriangleleft h),S)
y \tau_M.
\end{align*}
Therefore
\begin{align*} m(\alpha\beta)&=\tau^{-1}_Pg\lv(R,(x\vartriangleright
h),(x;h))\lh(X,\rv(R,(x\vartriangleright
h),(x;h))(x\vartriangleleft
h),S)\\&\quad\phi\left(\rh(X,\rv(R,(x\vartriangleright
h),(x;h))(x\vartriangleleft h),S)\right) \phi\left(y \right)\tau_M\\
&\hspace{-0.3cm}=\tau^{-1}_Pg\rv(R,(x\vartriangleright
h),(x;h))\phi\left(\rv(R,(x\vartriangleright
h),(x;h))\right)\phi\left(x\vartriangleleft
h\right)\phi\left(y\right)\tau_M\\
&\hspace{-0.3cm}=\tau^{-1}_Pg(x\vartriangleright
h)\phi\left(x\vartriangleleft h\right)\phi\left(y\right)\tau_M
\\& \hspace{-0.3cm}=\tau^{-1}_Pg\phi\left(x\right)h\phi\left(y\right)
\tau_M=m(\alpha)m(\beta).\\
\end{align*}
The second equality by \eqref{c2}, the third by \eqref{c1} and the
fourth by \eqref{c3}.

\end{proof}

\section{Examples}\label{exam}

In this section we shall give examples of braided groupoid data.

\subsection{ Handy braided groupoids}\label{handy}

\

In this section we study braided groupoid datum with the following
properties:
\begin{align}\label{cas1} &\bullet \, VH=HV, \\
\label{cas2} &\bullet \, \gamma(\Pc) H= H\gamma(\Pc), \text{ and }\\
\label{cas3} &\bullet \, \gamma(\Pc) V= V \gamma(\Pc).
\end{align}
This class of braided groupoids is the simplest to deal with. A
braided groupoid $\Vc$  whose associated  braided groupoid datum
$(D,V,H,\gamma)$ satisfies equations \eqref{cas1}, \eqref{cas2},
\eqref{cas3} will be called {\it handy braided groupoid}.

\medbreak

Let  $F$ be a group, and $\vartriangleright,
\vartriangleleft:F\times F\to F$ a left (respect.  right) action
on the set $F$. Let $\Pc$ be a set together with an operation
$\Pc\times\Pc\to \Pc, \, (P,Q)\mapsto PQ,$ not necessarily
associative, such that
\begin{itemize}
     \item[(i)]  There exists $O\in\Pc$ satisfying $PO=OP=P$, for all
    $P\in\Pc$,
   \item[(ii)] for any  $P\in\Pc$ there is a unique $Q\in\Pc$ such that
    $PQ=QP=O$. This element will be denoted by $P^{-1}$.

\end{itemize}

\medbreak

Let $\rightharpoondown:F\times \Pc\to \Pc$ be a group action and
$\s:\Pc\times\Pc\to F$ a map such that
\begin{align}\label{sg1} & \s(P,O)=\s(O,P)=1,\\
\label{sg2} & g\rightharpoondown O =O,\\
\label{sg3} & \s(P,P^{-1})=1,
\end{align}
for all $g\in F$, $P\in\Pc$.

\medbreak

Denote by $F\bowtie_{\s}\!\Pc_{\s\!}\bowtie F$ the  set
$F\times\Pc\times F$ with  multiplication given by
$$(g,P,x)(h,Q,y):=(g(x\vartriangleright h) \s(X,Y),XY,\s(X,Y)^{-1}
(x\vartriangleleft h)y),
$$ for every $g,h, x,y \in F, P,Q\in\Pc$, where $$X=
(x\vartriangleright h)^{-1}\rightharpoondown  P,\quad Y=
(x\vartriangleleft h)\rightharpoondown Q.$$

\medbreak

Equations \eqref{sg1}, \eqref{sg2} implies that $(1,O,1)$ is a
unit for this product.

\medbreak

Under certain compatibilities of the maps $\s, \triangleright,
\triangleleft, \rightharpoondown $ this multiplication ma\-kes
$F\bowtie_{\s}\!\Pc_{\s\!}\bowtie F$ into a group. This is the
next lemma.

\begin{lema}\label{tek} Keep the notation above. The set
$F\bowtie_{\s}\!\Pc_{\s\!}\bowtie F$ is a group with unit
$(1,O,1)$ if and only if the following conditions are fulfilled.
\begin{align}\label{dg1}&(F,F,\vartriangleright,
\vartriangleleft)\,\, \text{ is a matched pair of groups,}\\
\label{dg2}&(PQ)(\s(P,Q)^{-1}\rightharpoondown R)=(\s(Q,R)^{-1}
\rightharpoondown P)(QR),\\
\label{dg4}&\s(Q,R)\,\s\!\left(\s(Q,R)^{-1}\rightharpoondown
P,QR\right)=
\s(P,Q)\,\s\!\left(PQ,\s(P,Q)^{-1}\rightharpoondown R\right)\\
\label{dg6}&(g\rightharpoondown P)(g\rightharpoondown
Q)=\left(g\triangleleft \s(P,Q)\right) \rightharpoondown PQ,\\
\label{dg7}&g\triangleright
\s(P,Q)=\s(g\rightharpoondown P, g\rightharpoondown Q),\\
\label{dg8}&\left(g\triangleright
\s(P,Q)\right)\left(g\triangleleft \s(P,Q)\right)=g\,\s(P,Q),
\end{align}
for all $g\in F$, $P,Q, R\in\Pc$
\end{lema}

\begin{proof} Assume that
$F\bowtie_{\s}\!\Pc_{\s\!}\bowtie F$ is a group. From equalities
\begin{align*} (1,O,x)\left((1,O,y)(g,O,1)\right)&=\left((1,O,x)(1,O,y)\right)(g,O,1),\\
(1,O,x)\left((g,O,1)(h,O,1)\right)&=\left((1,O,x)(g,O,1)\right)(h,O,1),
\end{align*}
follow that $(F,F,\vartriangleright, \vartriangleleft)$ is a
matched pair of groups.  Equations \eqref{dg2}, \eqref{dg4} follow
from the equation
$$(1,P,1)\left((1,Q,1)(1,R,1)\right)=\left((1,P,1)(1,Q,1)\right)(1,R,1).$$
 Equations \eqref{dg6}, \eqref{dg7}, \eqref{dg8} can be
deduced from the equality
\begin{align*} (1,O,g)\left((1,P,1)(1,Q,1)\right)&=\left((1,O,g)(1,P,1)\right)(1,Q,1).\\
\end{align*}

\medbreak

Assume that equations \eqref{dg1} to \eqref{dg8} are fulfilled.
First we shall prove that the product in
$F\bowtie_{\s}\!\Pc_{\s\!}\bowtie F$ is associative. We claim that
it is enough to prove that
\begin{align}\label{as1}\left((g,P,1)(1,O,x) \right)(h,q,y)&=(g,P,1)\left((1,O,x)(h,q,y)
\right),\\
\label{as2}\left((1,O,x)(h,Q,y)\right)(f,R,z)&=(1,O,x)\left((h,Q,y)(f,R,z)\right),\\
\label{as3}\left((g,P,1)(h,Q,y)\right)(f,R,z)&=(g,P,1)\left((h,Q,y)(f,R,z)\right),
\end{align}
for al $P,Q,R\in\Pc$, $x,y,z,h,f,g\in F$. Indeed, let
$P,Q,R\in\Pc$, $x,y,z,h,f,g\in F$ then
\begin{align*} (g,P,x)\left((h,Q,y)(f,R,z)\right)&=\left((g,P,1)(1,O,x)\right)
\left((h,Q,y)(f,R,z)\right)\\
&=(g,P,1)\left((1,O,x)((h,Q,y)(f,R,z))\right)\\
&=(g,P,1)\left(((1,O,x)(h,Q,y))(f,R,z)\right)\\
&=\left((g,P,1)((1,O,x)(h,Q,y))\right)(f,R,z)\\
&=\left(((g,P,1)(1,O,x))(h,Q,y)\right)(f,R,z)\\
&=\left((g,P,x)(h,Q,y)\right)(f,R,z).
\end{align*}

The second equality by \eqref{as1}, the third by \eqref{as2}, the
fourth by \eqref{as3} and the fifth again by \eqref{as1}.

\medbreak

Equation \eqref{as1} follows by a direct calculation. Equation
\eqref{as2} follows from \eqref{dg1}, \eqref{dg7} and \eqref{dg8}.
Equation \eqref{as3} follows from \eqref{dg2},\eqref{dg4}
\eqref{dg6} and \eqref{dg7}. The inverse of an element is
$$(g,P,x)^{-1}=(x^{-1}\triangleright g^{-1}, (x^{-1}\triangleleft
g ^{-1})g\rightharpoondown P^{-1}), x^{-1}\triangleleft g^{-1}).
$$\end{proof}
When the map $\s$ or the action $\rightharpoondown$ are trivial,
conditions in Lemma \ref{tek} are easy to handle, as the following
corollaries show.

\begin{cor}\label{coro1} Assume that $(F,\triangleright,\triangleleft)$
is a matched pair of groups, $\Pc$ is a group with identity $O$,
and $\s:\Pc\times \Pc\to F$ is a map such that
\begin{align*} \bullet\quad &\s(P,O)=\s(O,P)=1,\\
\bullet\quad & \s(P,P^{-1})=1\\
\bullet\quad & \s(Q,R)\s(P,QR)= \s(P,Q)\s(PQ, R),
\end{align*}
for all $P,Q, R\in \Pc$. If in addition we have that
\begin{align*}
\bullet\quad & g\triangleright \s(P,Q)=\s(P,Q),\quad
g\triangleleft \s(P,Q)=\s(P,Q)^{-1}g\s(P,Q),
\end{align*}
for all $g\in V$, $P,Q\in\Pc$, then
$F_{\s\!}\bowtie\Pc_{\s\!}\bowtie F $ is a group, where
$\rightharpoondown$ is trivial.\qed
\end{cor}

\begin{cor}\label{coro2} Assume that $(F,\triangleright,\triangleleft)$
is a matched pair of groups, $\Pc$ is a group with identity $O$
and $\rightharpoondown$ is a left action of $F$ on $\Pc$ by group
automorphisms. Then $F\bowtie\Pc\bowtie F$ is a group, here the
map $\s$ is assumed to be trivial.\qed
\end{cor}

\medbreak

Let us assume that $F\bowtie_{\s}\!\Pc_{\s\!}\bowtie F$ is a
group, or, equivalently, the properties \eqref{dg1}, \eqref{dg2},
\eqref{dg4}, \eqref{dg6}, \eqref{dg7}, \eqref{dg8} hold.

\medbreak

Define the subgroups $V, H$ of $F\bowtie_{\s}\!\Pc_{\s\!}\bowtie
F$ by $V:=F\times O\times 1,\, H:=1\times O\times F$. The map
$\gamma:\Pc\to F\bowtie_{\s}\!\Pc_{\s\!}\bowtie F$, is the
inclusion; $\gamma(P)=(1,P,1)$.

Then the collection $(F\bowtie_{\s}\!\Pc_{\s\!}\bowtie F,V,
H,\gamma)$ satisfies conditions of Lemma \ref{l2} (iii).

\begin{teo}\label{mainex} If $(F,\vartriangleright, \vartriangleleft)$ is a
braided group, then $(F\bowtie_{\s}\!\Pc_{\s\!}\bowtie F,V,
H,\gamma)$ is a braided groupoid datum and the associated braided
groupoid is handy.

Reciprocally if $(D,V,H,\gamma,\phi)$ is a braided groupoid datum
and the associated braided groupoid is handy, then
$(V,\vartriangleright, \vartriangleleft)$ is a braided group,
$\Pc$ has an operation that satisfies (i), (ii), there are maps
$\s:\Pc\times\Pc\to V$, $\rightharpoondown: V\times \Pc\to\Pc$
such that $D\cong V\bowtie_{\s}\!\Pc_{\s\!}\bowtie V$ and $\gamma$
is the inclusion via this isomorphism.

\end{teo}

\begin{proof} If $h,y\in F$, $P,Q\in \Pc$ then
\begin{align*} \lambda_V(P,(h,1,1),Q)&=(h\s(h^{-1}\rightharpoondown P,
Q),1,1),\\
\rho_V(P,(h,1,1),Q)&=(1,1,\s(h^{-1}\rightharpoondown P,
Q)^{-1}),\\
\lambda_H(P,(1,1,y),Q)&=(\s(P,y\rightharpoondown Q),O,1),\\
\rho_H(P,(1,1,y),Q)&=(1,O,\s(P,y\rightharpoondown Q)^{-1} y),\\
(1,1,y)\triangleright (h,1,1)&= (y\triangleright h,1,1),\\
(1,1,y)\triangleleft (h,1,1)&=(1,1,y\triangleleft h).
\end{align*}
Therefore the first assertion follows from Theorem \ref{main}.
\medbreak

Let $(D,V,H,\gamma)$ be a braided groupoid datum such that
equations \eqref{cas1}, \eqref{cas2}, \eqref{cas3} are satisfied.
Abusing of the notation we define
$\triangleright,\triangleleft:V\times V\to V$ by
$$g\triangleright h:=\phi^{-1}(g)\triangleright h, \quad g\triangleleft h:=
\phi\left(\phi^{-1}(g)\triangleleft h\right), $$ for all $g,h\in
V$. Since $VH=HV$ then $(x;g)=O$ for all $x\in H, g\in V$.
Associativity axiom of the group $D$ implies that
$(V,V,\triangleright,\triangleleft)$ is a matched pair of groups.
Equation \eqref{c3} implies that $(V,\vartriangleright,
\vartriangleleft)$ is a braided group.

\medbreak

Define the following composition $\Pc\times\Pc\to\Pc$,
$PQ:=(P;1;Q)$. Clearly $O$ is a unit for this operation. The
existence of inverse in $D$ translates in the existence of the
inverse in $\Pc$.

\medbreak

Define the maps $\s:\Pc\times\Pc\to V$, $\rightharpoondown:
V\times \Pc\to \Pc$ by
 $$\s(P,Q):=\lv(P,1,Q),\quad  g\, \gamma_P:=\gamma_{g\rightharpoondown P} g'  $$
for all $P, Q\in \Pc$, $g\in V$, where $g'$ is some element in $G$
that depends on $g$ and $P$. Since the map $m$ is a groupoid
morphism then $m( g\, \gamma_P)=m(g)=g$, and therefore $g=g'$.
Hence the map $\rightharpoondown$ is defined by the equation
$$ g\,\gamma_P :=\gamma_{g\rightharpoondown P} \, g.$$

Equation \eqref{c1} implies that
$\rv(P,1,Q)=\phi^{-1}(\s(P,Q)^{-1})$.

\medbreak

Define $f:D\to V\bowtie_{\s\!}\Pc_{\s\!}\bowtie V$ by
$$f(g\gamma_P x)=(g,P,\phi(x)), $$
for all $g\in V, P\in\Pc, x\in H$. This is a well defined group
isomorphism. This ends the proof of the theorem.
\end{proof}

\medbreak

In particular, Theorem \ref{mainex}, in presence of Corollaries
\ref{coro1}, \ref{coro2}, shows that there is a way to produce
many examples of braided groupoid datum. For example, take
$(F,\triangleright,\triangleleft)$ any braided group, $\Pc$ a
group such that $F$ acts on $\Pc$ by group automorphism; or take
$F, \Pc$ two groups with a normalized 2-cocycle $\s:F\times F\to
\Pc$,  $\triangleright:F\times F\to F$ the trivial action and
$\triangleleft:F\times F\to F$ the adjoint action.

\begin{cor} Let $(D,V,H,\gamma)$ be a braided groupoid datum, where $V$
and $H$ are normal subgroups of $D$. Then the associated braided
groupoid $\Vc$ is handy, moreover the action $\rightharpoondown$
is trivial.
\end{cor}
\begin{proof} Since $V$ is normal
$\gamma_P g\gamma^{-1}_P\in V, $ for all $P\in\Pc,$ $g\in V$.
Hence, $\gamma_P g=g^{'}\gamma_P$ for some $g^{'}\in V$. Since
$(D,V,H,\gamma)$ is a braided groupoid datum then $g=g^{'}$.
Analogously we prove that $\gamma_P x=x\gamma_P$ and $gx=xg$ for
all $x\in H$, $g\in V$, $P\in \Pc$ .
\end{proof}

\subsection{Non-handy examples }\label{ej2}
Let $(A, \triangleright,\triangleleft)$ be a matched pair of
groups. Let $ \Pc$ be a group, and let $\psi:A\times A\to
\cent(\Pc)$, $\cent(\Pc)$ the center of $\Pc$, be a map such that
for any $a, b, c \in A$
\begin{align}\label{d1} \psi(a,bc)&=\psi(a,b)\psi(a\triangleleft b,c),\\
\label{d2} \psi(ab,c)&=\psi(a,b\triangleright c)\psi(b,c).
\end{align}

Define the group $D$ whose underlying set is $A\times \Pc\times A$
 and  multiplication given by
$$(a,P,c)(x,Q,z)=(a(c\triangleright x),P\psi(c,x)Q,(c\triangleleft x)z), $$
for any $a,c,x,z \in A,\,  P,Q \in \Pc$. A straightforward
computation shows that this operation is associative.

\medbreak

Let $V=A\times 1\times 1,$ $H=1\times 1\times A$ and
$\gamma:\Pc\to D$, $\gamma_P=(1,P,1)$,

\begin{lema} If $(A, \triangleright,\triangleleft)$ is a braided
group then the collection $(D,V,H,\gamma)$ is a braided groupoid
datum.
\end{lema}
\begin{proof} For any $P, Q \in \Pc$, $a,b\in A$ we have that
\begin{align*} &\lambda_V(P,(a,1,1),Q)=(a,1,1),\quad
\rho_V(P,(a,1,1),Q)=1\\
&\lambda_H(P,(1,1,a),Q)=1,\quad \rho_H(P,(1,1,a),Q)=(1,1,a),\\
&(1,1,a)\triangleright (b,1,1)=(a\triangleright b,1,1),\quad
(1,1,a)\triangleleft (b,1,1)=(1,1,a\triangleleft b).
\end{align*}
Then the Lemma follows by applying Theorem \ref{main}.
\end{proof}

If $a,z\in A$ then
$$ (a,1,1)(1,1,z)= (a,1,z), \quad  (1,1,z)(a,1,1)=(z\triangleright a,
\psi(z,a),z\triangleleft a). $$ Thus, $VH=HV$ if and only if
$\psi=1$.

\begin{rmk} There are many collections $(A, \triangleright,
\triangleleft, \psi)$, where $(A, \triangleright, \triangleleft)$
is a braided group and $\psi$ is a map satisfying \eqref{d1},
\eqref{d2}. For example take $A$ any group, $\triangleright$ the
adjoint action, $\triangleleft$ the trivial action and $\psi$ any
bicharacter, that is $\psi:A\times A\to \cent(\Pc)$ such that
\begin{align*} \psi(a,bc)&=\psi(a,b)\psi(a,c),\\
 \psi(ab,c)&=\psi(a,c)\psi(b,c),
\end{align*}
for all $a,c,x,z \in A$.
\end{rmk}

This class of examples arise from the following general
observation.  Let $(D,V,H,\gamma,\phi)$ be any braided groupoid
datum. Recall the map $(\, ;\, ):V\times H\to \Pc$ defined by
equation \eqref{b1}. If we assume that for all $P\in\Pc, g\in G,
x\in H$
$$\gamma_P g= g \gamma_P,\quad \gamma_Px=x\gamma_P,$$
then the map $\psi:V\times V\to \Pc$ defined by
$$\psi(g,h)=(g;\phi^{-1}(h)), $$
for all $g,h\in V$, satisfies equations  \eqref{d1} and
\eqref{d2}. Consider the following operation in $\Pc$;
$P.Q=(P;1;Q)$. Since $xg\gamma_P= \gamma_Pxg$ for all $x\in H,g\in
V, P\in \Pc$ then
$$(x\triangleright g)\gamma_{(x; g)}\gamma_P (x\triangleleft g)=
(x\triangleright g)\gamma_P\gamma_{(x; g)} (x\triangleleft g),$$
 and thus, $(x; g)\in \cent(\Pc)$ for all $x\in H,g\in
V$.

\section{The Braiding}\label{section:trenza}

In this section we explicitly compute the braiding for the braided
data given in the previous section.

\medbreak

Let $(D,V,H,\gamma,\phi)$ be a braided groupoid datum and let
$\D=\Vc\Hc$ be the associated exact factorization of groupoids.
Let $\alpha\in\Hc, \beta \in\Vc$ then
$$\alpha=\tau^{-1}_P\, \gamma_P x
\gamma^{-1}_Q\,\tau_Q,\quad \beta=\tau^{-1}_Q\,g \,\tau_R, $$ for
some $P,Q,R\in\Pc, g\in V, x\in H$. Then
$$\alpha\beta=\tau^{-1}_P\, \gamma_P x
\gamma^{-1}_Qg \,\tau_R.$$ Since
$\alpha\beta=(\alpha\rightharpoonup\beta)(\alpha\leftharpoonup\beta)$,
the determination of the actions $\rightharpoonup, \leftharpoonup$
relies on the explicit calculation of $\gamma_P x \gamma^{-1}_Qg$.
This will be done in the following for the examples explained
above.

\subsection{The braiding for handy braided groupoids}\label{calculo}

Let $\Vc$ be a handy braided groupoid and
$(F\bowtie_{\s\!}\Pc_{\s\!}\bowtie F,V, H,\gamma)$ its braided
groupoid datum. Let also $\D=\Vc\Hc$ be the exact factorization
associated to the collection $(F\bowtie_{\s\!}\Pc_{\s\!}\bowtie
F,V, H,\gamma)$.

\begin{lema}\label{te3} IF $P,Q, R\in\Pc, x,y\in F$ then
\begin{align*} \gamma_P(1,O,x)\gamma^{-1}_Q (y,O,1)\gamma_R=&
(\s(P,Q^{-1})(\s(P,Q^{-1})^{-1}x\triangleright y)\s(S,T),O,1)\gamma_{ST}\\
&(1,O,\s(P,Q^{-1})^{-1}x\triangleleft y),
\end{align*}
where
\begin{align}\label{s} S&=(\s(P,Q^{-1})^{-1}\triangleright y)^{-1}\rightharpoondown
(PQ^{-1}),\\
\label{t} T&=(\s(P,Q^{-1})^{-1}\triangleright y)\rightharpoondown
R.
\end{align}
\end{lema}
\begin{proof} Straightforward.
\end{proof}

\medbreak

Let $(\alpha, \beta)\in\Vc_{\,\tgt\!}\times_{\sou}\Vc$. Then there
exists $P,Q,R\in \Pc$, $x,y\in F$ such that
$$\alpha=\tau^{-1}_P(x,O,1)\tau_Q,\quad \beta=\tau^{-1}_Q (y,O,1)\tau_R. $$
Then
\begin{align*} \phi^{-1}(\alpha)\beta&=\tau^{-1}_P\gamma_P
(1,O,x)\gamma^{-1}_Q  (y,O,1) \tau_R\\
&=\tau^{-1}_P\gamma_P (1,O,x)\gamma^{-1}_Q
(y,O,1)\gamma_R\gamma^{-1}_R \tau_R\\
&=\tau^{-1}_P (\s(P,Q^{-1})(\s(P,Q^{-1})^{-1}x\triangleright
y)\s(S,T),O,1) \tau_{ST} \\
&\quad\,\,
\tau^{-1}_{ST}\gamma_{ST}(1,O,\s(P,Q^{-1})^{-1}x\triangleleft
y)\iota^{-1}_R \tau_R,
\end{align*}
Where $S, T\in\Pc$ are as in Lemma \ref{te3}. Since
$$\tau^{-1}_P \,(\s(P,Q^{-1})(\s(P,Q^{-1})^{-1}x\triangleright
y)\s(S,T),O,1) \,\tau_{ST}\in \Vc(P,ST),
$$
$$\tau^{-1}_{ST}\,\gamma_{ST}(1,O,\s(P,Q^{-1})^{-1}x\triangleleft
y)\gamma^{-1}_R\, \tau_R\in\Hc(ST,R), $$ then

\begin{align*} \alpha\rightharpoonup \beta &=\tau^{-1}_P\,
(\s(P,Q^{-1})(\s(P,Q^{-1})^{-1}x\triangleright
y)\s(S,T),O,1)\, \tau_{ST}\\
\alpha\leftharpoonup\beta
&=\tau^{-1}_{ST}\,(1,O,\s(P,Q^{-1})^{-1}x\triangleleft y)\, \tau_R
\end{align*}
As a consequence of these calculations we have the following
result.

\begin{prop}\label{ss} The braiding for the handy braided groupoid $\Vc$ is
given by the formula
\begin{align*} c(\alpha,\beta)=&(\tau^{-1}_P\,
(\s(P,Q^{-1})(\s(P,Q^{-1})^{-1}x\triangleright y)\s(S,T),O,1)\,
\tau_{ST},\\ &\tau^{-1}_{ST}\,(1,O,\s(P,Q^{-1})^{-1}x\triangleleft
y)\, \tau_R)
\end{align*}
where $\alpha=\tau^{-1}_P(x,O,1)\tau_Q, \beta=\tau^{-1}_Q
(y,O,1)\tau_R$ and $S, T$ are given by equations \eqref{s},
\eqref{t} \qed.
\end{prop}

\begin{rmk} When $\# \Pc=1$ then formula in Proposition \ref{ss} is
$c(x,y)=(x\triangleright y, x\triangleleft y)$, which is the braid
formula for the braided group $(F,\triangleright,\triangleleft)$.
\end{rmk}

\subsection{The braiding for the examples in subsection \ref{ej2} }

Let $(A, \triangleright ,\triangleleft)$ be a braided group, $\Pc$
be a group. Let also $\psi:A\times A\to \cent(\Pc)$ be a map
satisfying \eqref{d1}, \eqref{d2}. Let $(D,V,H,\gamma)$ be the
braided groupoid datum as in example \ref{ej2}. Let $\D=\Vc\Hc$ be
the exact factorization associated to $(D,V,H,\gamma)$.

\begin{lema} Let $a, b\in A$, $P,Q, R\in\Pc$ then
$$\gamma_P(1,1,a)\gamma^{-1}_Q(b,1,1)\gamma_R=(a\triangleright b,
\psi(a,b)PQ^{-1}R,a\triangleleft b). \qed$$
\end{lema}

Let $(\alpha, \beta)\in\Vc_{\,\tgt\!}\times_{\sou}\Vc$. Then there
exists $P,Q,R\in \Pc$, $a,b\in A$ such that
$$\alpha=\tau^{-1}_P (a,1,1)\tau_Q,\quad \beta=\tau^{-1}_Q(b,1,1)\tau_R, $$
Then
\begin{align*} \phi^{-1}(\alpha)\beta&=\tau^{-1}_P\gamma_P(1,1,a)
\gamma^{-1}_Q(b,1,1)\tau_R\\
&=\tau^{-1}_P\gamma_P(1,1,a) \gamma^{-1}_Q(b,1,1)\gamma_R\gamma^{-1}_R\tau_R\\
&=\tau^{-1}_P \,(a\triangleright
b,\psi(a,b)PQ^{-1}R,a\triangleleft b)
\gamma^{-1}_R\,\tau_R\\
&=\tau^{-1}_P\, (a\triangleright b,1,1) \,\tau_S\, \, \tau^{-1}_S
\gamma_S (1,1,a\triangleleft b)\gamma^{-1}_R\,\tau_R,
\end{align*}
where $S=\psi(a,b)PQ^{-1}R\in\Pc$. Since
$$\tau^{-1}_P\, (a\triangleright b,1,1)\, \tau_S\in \Vc(P,S), $$
$$ \tau^{-1}_S \gamma_S
(1,1,a\triangleleft b)\gamma^{-1}_R\,\tau_R\in \Hc(S,R),$$ then

\begin{align*} \alpha\rightharpoonup\beta&=\tau^{-1}_P\, (a\triangleright b,1,1)\,
\tau_S\\
\alpha\leftharpoonup\beta&=\tau^{-1}_S (1,1,a\triangleleft
b)\,\tau_R.
\end{align*}

\begin{prop} If $\alpha=\tau^{-1}_P
(a,1,1)\tau_Q,\, \beta=\tau^{-1}_Q(b,1,1)\tau_R$ then the braiding
for the examples \ref{ej2} are given by the formula
$$c(\alpha, \beta)=(\tau^{-1}_P\, (a\triangleright b,1,1)\,
\tau_S, \tau^{-1}_S (1,1,a\triangleleft b)\,\tau_R), $$ where
$S=\psi(a,b)PQ^{-1}R$.\qed
\end{prop}


\begin{thebibliography}{AJS}

\bibitem[A]{A}  {\sc N. Andruskiewitsch},
{\it On the quiver-theoretical quantum Yang-Baxter equation},
Selecta Math.(N.S.) to appear  \texttt{math.QA/0402269}.

\bibitem[AA]{AA}  {\sc M. Aguiar} and {\sc N. Andruskiewitsch},
{\it Representations of matched pairs of groupoids and
applications to weak Hopf algebras}, preprint, (2004),
\texttt{math.QA/0402118}, Contemp. Math. to appear.

\bibitem[AM]{AM}  {\sc N. Andruskiewitsch } and {\sc J.M. Mombelli},
{\it Examples of weak Hopf algebras arising from vacant double
groupoids}, \texttt{math.QA/0405374}, submitted.


\bibitem[AN]{AN}  {\sc N. Andruskiewitsch } and {\sc S. Natale},
{\it Double categories and quantum groupoids},
\texttt{math.QA/0308228}, Publ. Mat. Uruguay, to appear.

\bibitem[B]{B}  {\sc E. Beggs}, {\it Making non-trivially
associated tensor categories from left coset representatives}, J.
Pure Appl. Algebra {\bf 177} (2003), 5--41.


\bibitem[D]{D}  {\sc V.G. Drinfeld},
{\it On some unsolved problems in quantum group theory},  Lect.
Notes  Math. {\bf 1510}, Springer-Verlag, Berlin (1992).

\bibitem[EGS]{EGS}  {\sc P. Etingof}, {\sc R. Guralnik} and {\sc A.
Soloviev}, {\it Indecomposable set-theoretical solutions to the
quantum Yang-Baxter equation on a set with prime number of
elements}, J. Algebra {\bf 242} 2 (2001), 709--719.


\bibitem[ESS]{ESS}  {\sc P. Etingof}, {\sc T. Schedler} and {\sc A.
Soloviev}, {\it Set-theoretical solutions to the quantum
Yang-Baxter equation}, Duke Math. J. {\bf 100} (1999), 169--209.

\bibitem[LYZ1]{LYZ} {\sc Jiang-Hua Lu}, {\sc Min Ya} and {\sc
Yong-Chang Zhu}, {\it On the set-theoretical Yang-Baxter
equation}, Duke Math. J. {\bf 104} (2000), 1--18.

\bibitem[LYZ2]{lyz3} {\sc Jiang-Hua Lu, Min Yan} and {\sc Yong-Chang Zhu},
\emph{Quasi-triangular structures on Hopf algebras with positive
bases\/}, in ``New trends in Hopf Algebra Theory"; Contemp. Math.
\textbf{267} (2000), 339--356.


\bibitem[Ma]{Ma} {\sc K. Mackenzie},  \emph{Double Lie algebroides
and second-order geometry I}, Adv. Math. {\bf 94} (1992), pp.
180--239.


\bibitem[N]{N} {\sc S. Natale}, \emph{Frobenius-Schur indicators for a
class of fusion categories}, to appear in Pacific J. Math.
Preprint math.QA/0312466.

\bibitem[S]{S} {\sc A.Soloviev}, \emph{Non-unitary set-theoretical
solutions to the quantum Yang-Baxter equation}, Math. Res. Lett.
{\bf 7} (2000), no. 5-6, pp.577--596.

%\bibitem[R]{R} {\sc J. Renault}, {\it A groupoid approach to
%$C^{*}$-algebras}, Lect. Notes  Math. {\bf 793}, Springer-Verlag,
%Berlin (1980).

\bibitem[T]{T}  {\sc M. Takeuchi, }
{\it Survey on matched pairs of groups. An elementary approach to
the ESS-LYZ theory,}  Banach Center Publ. {\bf 61} (2003),
305--331.


\end{thebibliography}
\end{document}